# EXISTENCE OF PHI-RELATIONSHIP WITHIN REGULAR HEXAGONAL TESSELLATION

E. Bervalds

Latvian Academy of Sciences
1 Akademijas laukums, Riga, LV-1050, LATVIA
e-mail: berv@acad.latnet.lv

The Phi- relationship also known as Phi-factor appears in a number of lattice structures, mostly considering the lines within several separate circles or polygons. The paper considers a regular hexagonal tessellation as a lattice with the highest specific mechanical stiffness."

## 1. THE Phi CENTER OF A THREE-ELEMENT CLUSTER OF HEXAGONAL TESSELLATION

Any straight line segment can be divided into two parts so that the ratio between the total length of the parts and the longest part is equal to the ratio between the longest part and the shortest one. This is possible provided that the division point of the line segment is concurrent with the Phi center which makes the values of both the ratios equal to 1.618…[1].

We can now consider a separate cluster of hexagonal {6;3} tessellation consisting of three elements with common point O (Fig. 1). Each hexagon will be encircled with three concentric circles with the radii relationship of 1:2:4, the middle circle being the encircling one for each separate hexagon. If we draw six tangents of the three smallest circles from point O to their crossing points with the middle radii circles we will obtain six straight line segments of equal length $A_k B_k (k = 1-6)$ having the point O as the Phi relationship division point [2]. Simultaneously, at this construction point O there is not only a Phi center for any of $k$ line segments, but also the Phi center of a selected hexagonal cluster supported by equal values of the six separate ratios, i.e.

$$\frac{A_k B_k}{A_k O} = \frac{A_k O}{O B_k} = 1{,}6180339887, \qquad (1)$$

Similar constructive processes can be repeated by choosing any vertex of any hexagon in the tessellation as the center of a separate cluster. All these centers could be defined as the Phi centers of the hexagonal tessellation according to previously described Phi relationships in the six directions.

The level of the Phi relationship determined in (1) achieved within a regular hexagonal tessellation is assessed by placing it in the convergence sequence of Fibonacci's numbers
on Phi [3]. This corresponds to the 11-th series with the ratio 1,6181818181 by the 0,0001478294 variance from Phi.

## 2. CONCLUSION

Any separate three-element cluster of a regular hexagonal tessellation sets the vertices of a hexagon as the Phi centers with the variance from Phi that corresponds to the 11-th series of the Fibonacci convergence sequence on Phi.

*Fī* ATTIECĪBU ESAMĪBA
REGULĀRU SEŠSTŪRU MOZAĪKĀ

E. Bervalds

K o p s a v i l k u m s

Ģeometriskas konstruēšanas ceļā pierādīta seškārtīga *Fī* attiecību esamība atsevišķā trīs sešstūru kopā, definējot to kopējo virsotni kā *Fī* centru un vispārinot šo apzīmējumu uz jebkuru regulāru sešstūru mozaīkas virsotni. Iegūtās *Fī* skaitliskās vērtības precizitāte novērtēta to salīdzinot ar atbilstošo Fibonacci skaitļu rindas locekli, tai konverģējot uz *Fī*.
19.03.2007.

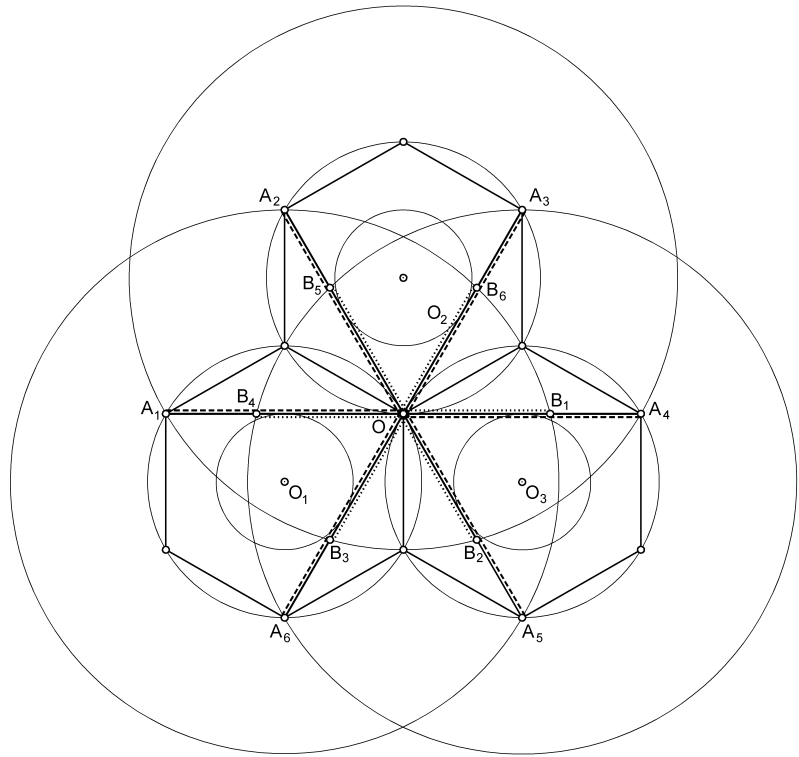

*Fig. 1.* Geometrical structure of three-element cluster of the hexagonal
tessellation for determination of the Phi center